\documentclass[12pt]{article}
\usepackage{graphicx}
\usepackage{amsfonts, amsmath,  amssymb, latexsym}
\usepackage{eucal}

\newcommand\comment[1]{}                

\def\int{\mathop{\rm int}}

\newtheorem{theorem}{Theorem}

\newtheorem{lemma}{Lemma}[section]

\newenvironment{definition}[1]{\par \medskip \par \noindent
\texttt{ Definition:}\ #1}{ \par \bigskip}

\renewcommand{\phi}{\varphi}          
\renewcommand{\epsilon}{\varepsilon}
\renewcommand\tilde{\widetilde}

\newcommand\conv{\operatorname{conv}}
\newcommand\aff{\operatorname{aff}}
\newcommand\lin{\operatorname{lin}}

\renewcommand\vert{\operatorname{vert}}

\newcommand\Ker{\operatorname{Ker}}
\newcommand\GCD{\operatorname{GCD}}

\newcommand\uu{\mathbf{u}}  
\newcommand\vv{\mathbf{v}} 
\newcommand\cc{\mathbf{c}}

\newcommand\x{\mathbf{x}}

\newcommand\cP{\mathcal{P}}

\newcommand\cE{\mathcal{E}}

\newcommand\R{\mathbb{R}}    
\newcommand\N{\mathbb{N}}    
\newcommand\Z{\mathbb{Z}}    

\input{tcilatex}

\begin{document}

\author{Robert Erdahl, Andrei Ordine, and Konstantin Rybnikov}
\title{Perfect Delaunay polytopes and Perfect Inhomogeneous Forms}
\date{June 11, 2004 }
\maketitle

\begin{abstract}
A lattice Delaunay polytope $D$ is called \textit{perfect} if it has the
property that there is a unique circumscribing ellipsoid with interior
free of lattice points, and with the surface containing only those lattice
points that are the vertices of $D$. \ An inhomogeneous quadratic form is
called \textit{perfect} if it is determined by such a circumscribing
''empty ellipsoid'' uniquely up to a scale factor. Perfect inhomogeneous
forms are associated with perfect Delaunay polytopes in much the way that
perfect homogeneous forms are associated with perfect point lattices. \ We
have been able to construct some infinite sequences of perfect Delaunay
polytopes, one perfect polytope in each successive dimension starting at
some initial dimension; we have been able to construct an infinite number
of such infinite sequences. \ Perfect Delaunay polytopes are intimately
related to the theory of Delaunay polytopes, and to Voronoi's theory of
lattice types.
\end{abstract}

\section{ \ Introduction\label{introduction}}

Consider the lattice $\mathbb{Z}^{d}$, and a lattice polytope $D$.  If $D$
can be circumscribed by an ellipsoid $\mathcal{E}=\partial \{\,\x \in \R^d
\; \vline \; f(\x) \le R^2   \}$, where $f$ is a quadratic function in
$x_1,\ldots,x_d$, with no
 $\mathbb{Z}^{d}$-elements interior to $\{\,\x \in \R^d
\;\vert\; f(\x) \le R^2   \}$, and with all $\mathbb{Z}^{d}$-elements on
$\cE$ being the vertices of $D$, we will say that $D$ is a
\textit{Delaunay} polytope with respect to the metric form defined by
$\cE$; more informally, we will say that $D$ is \textit{Delaunay} if it
can be circumscribed by an \textit{empty ellipsoid} $\mathcal{E}$.
Typically there is a family of empty ellipsoids that can be circumscribed
about a given Delaunay polytope $D$ but, if there is only one, so that
$\mathcal{E}$ is uniquely determined by $D$, we will say that $D$ is a
\textit{perfect Delaunay polytope}. Perfect Delaunay polytopes are
distinguished in that they cannot fit properly inside other lattice
Delaunay polytopes. Perfect Delaunay polytopes are fascinating geometrical
objects -- examples are the six and seven dimensional Gossett polytopes
with $27$ and $56$ vertices known as $2_{21}$ and $3_{21}$ in Coxeter's
notation, which appear in the Delaunay tilings for the quadratic forms
corresponding to the root lattices $E_{6}$ and $E_{7}$.

We have studied perfect Delaunay polytopes by constructing infinite
sequences of them, one perfect Delaunay polytope in each successive
dimension, starting at some initial dimension.  We have been able to
construct an infinite number of infinite sequences of perfect Delaunay
polytopes. One of our constructions is the sequence of \textit{ G-topes},
G$^{d},d=6,7,...$ , with the initial term being the six-dimensional
Gossett polytope $2{21}$ with $27$ vertices; each G-tope is asymmetric
with respect to central inversion, and G$^{d}$ has $\binom{d+2}{2}-1$\
vertices. \ \ Another construction is the sequence of \textit{C-topes},
$\textrm{C}^{d},d=7,8,...$ , with initial term the seven-dimensional
Gossett polytope $3_{31}$ with $56$ vertices; each C-tope is symmetric
with respect to central inversion, and C$^{d}$ has $2 \binom{d+1}{2}$
vertices.\ Just as the six-dimensional Gossett polytope can be represented
as a section of the seven dimensional one, each term $\textrm{G}^{d}$ of
the asymmetric sequence can be represented as a section of the term
$\textrm{C}^{d+1}$ of the symmetric sequence. \  \

Each perfect Delaunay polytope in the sequence of G-topes uniquely
determines a lattice that is an analogue of the root lattice $E_{6}$, and
each term in the sequence of C-topes uniquely determines a lattice that is
an analogue of the root lattice $E_{7}$. \ \ The Voronoi and Delaunay
tilings for the lattices in these sequences of lattices show many features
of the lead terms, namely, the Voronoi and Delaunay tilings for the
quadratic forms corresponding to the root lattices $E_{6}$ and $E_{7}$.

There is a number of reasons why perfect Delaunay polytopes are
fascinating objects for study. First, as mentioned before, they are
"perfect" inhomogeneous analogs of perfect lattices. This alone seems to
be a natural motivation for a geometer of numbers.  If the empty ellipsoid
$\cE$ circumscribing a perfect Delaunay polytope $D$ with the vertex set
$\vert D$ is defined by an equation $f(\x)=0$, where $f$ is an
inhomogeneous quadratic form, then all the coefficients of $f$ are
uniquely determined (up to a common scaling factor) from the system
$\{f(\x)=1\: | \: \x \in \vert D\}$. Clearly, $\vert D$ is then the set of
all points on which $f$ attains its minimum over $\Z^d$. An analogous
property for homogeneous quadratic forms is called \emph{perfection} after
Korkin and Zolotareff (1873) \comment{\cite{KZ}} (see Martinet (2003) and
\comment{\cite{Mart,CSperf}} Conway and Sloane (1988) for modern treatment
of perfect homogeneous forms). Naturally, in our context $f$ is called an
\emph{inhomogeneous perfect form} and $\cE$ is referred to as a
\emph{perfect ellipsoid}. The vertices of a perfect Delaunay polytope are
analogs of the minimal vectors for a perfect form: the minimal possible
number of vertices of a perfect Delaunay polytope is $\frac{n(n+1)}{2}+n,$
while the minimal number of shortest vectors of a perfect forms is
$\frac{n(n+1)}{2}$.

Perfect Delaunay polytopes were first considered by Erdahl in 1975 in
connection with lattice polytopes arising from the quantum mechanics of
many electrons.  He showed (1975) that there are perfect Delaunay
polytopes in one-dimensional lattices, and showed that the Gosset polytope
G$^6=2_{21}$ with 27 vertices was a perfect Delaunay polytope in the root
lattice $E_{6}$. \ He also showed that there were no perfect Delaunay
polytopes in dimensions 2, 3, and 4. \ These results were extended by
Erdahl (1992) by showing that the 7-dimensional Gosset polytope
C$^7=3_{31}$ with 56 vertices is perfect, and that there are no perfect
Delaunay polytopes in lattices with dimension lying between one and six.
Erdahl also proved that $[0,1]$, G$^6$, and C$^7$ are the only perfect
Delaunay polytopes  existing in root lattices. Deza, Grishukhin, and
Laurent (1992,1995) found more examples of perfect Delaunay polytopes in
dimensions 15, 16, 22, and 23. The first construction of infinite
sequences of perfect Delaunay polytopes was first announced at the
Conference dedicated to the Seventieth Birthday of Sergei Ryshkov (Erdahl,
2001), and later reported by Rybnikov (2001) and Erdahl and Rybnikov
(2002). Perfect Delaunay polytopes have been classified up to dimension 7
-- Dutour (2004)  proved that G$^6=2_{21}$ is the only Gosset polytope for
$d=6$. It is strongly suspected that the existing lists of seven and eight
dimensional perfect Delaunay polytopes are complete (see
http://www.liga.ens.fr/$\tilde{~}$dutour/).

Perfect Delaunay polytopes are important in the study of \emph{strongly
regular graphs} or, more generally, \emph{theory of association schemes}.
In fact, the Schlafli graph and the Gosset graph can be constructed as
1-skeletons of Gosset polytopes  G$^6$ and C$^7$, and that is precisely
how the Gosset graph has been discovered. A \emph{maximal family of
equiangular lines} is a metrical concept, closely related to the
combinatorial notion of strongly regular graph (see Deza et al (1992),
Deza and Laurent (1997)  \comment{\cite{DGL,DL}} for details). All
Delaunay polytopes found by Deza et al \comment{(\cite{DGL})} have been
derived from such families of lines.

Perfect Delaunay polytopes play an important role in the \emph{ L-type
reduction theory of Voronoi and Delaunay} for positive (homogeneous)
quadratic forms. An L-type domain is the collection of all possible
quadratic metric forms that give the same Delaunay tiling $\mathcal{D}$
for $\mathbb{Z}^{d}$. L-type domains are relatively open polyhedral cones,
with boundary cells that are also L-type domains -- these conical cells
fit together to tile the cone of metrical quadratic forms, which is
described in the next section in  detail. Simplicial Delaunay tilings
label the full dimensional conic "tiles'', and all other possible Delaunay
tilings label the lower dimensional cones. \ A significant role for
perfect Delaunay polytopes is that they provide labels for a subclass of
\textit{ edge forms for L-types}. Edges-forms are forms lying on extreme
rays of full-dimensional L-type domains, which are labels by simplicial
Delaunay tilings. Prior to the discovery of the sequence
$\{\textrm{G}^{d}\}$ of G-topes by Erdahl and Rybnikov in 2001 only
finitely many edge forms were known. \ All of the infinite sequences of
perfect Delaunay polytopes that we report on correspond to infinite
sequences of edge forms for L-type domains. The significance of extreme
$L$-types is much due to their relation to the structure of Delaunay and
Voronoi tilings for lattices. The Delaunay tilings that correspond to edge
forms play an important role: \ \emph{The Delaunay tiling for an L-type
domain is the intersection of the Delaunay tilings for edge forms for the
L-type.} (Erdahl, 2000) \ There is a corresponding dual statement on the
structure of Voronoi polytopes: \emph{The Voronoi polytope }$V_{\varphi }$%
\emph{\ for a form }$\varphi $\emph{\ contained in an L-type domain is a
weighted Minkowski sum of linear transforms of Voronoi polytopes for each
of the edge forms. }\ \ The latter dual result was first established by
H.-F. Loesch  in  his 1990 doctoral dissertation, although it was first
published by Ryshkov (1998, 1999), who independently rediscovered Loesh's
theorem; this dual result was given a shorter and simpler proof by Erdahl
(2000). \

Edge forms that are interior to the cone of metric forms are rare in low
dimensions. They first occur in dimension 4: $D_{4}$ has an extreme $L
$-type. A good proportion, but not all, of the edge forms appearing in
lower dimensions relate either directly or indirectly to perfect Delaunay
polytopes. As shown by Dutour and Vallentin (2003) this situation does not
persist in higher dimensions: there is an explosion of edge forms in six
dimensions, and that only a tiny fraction of these are inherited from
perfect Delaunay polytopes.

 All of our sequences of perfect Delaunay polytopes determine
corresponding sequences of lattices with similar combinatorial properties.
Our constructions are first steps of a program to explore the geometry of
higher dimensional lattices through infinite sequences of lattices.  The
infinite sequences we have constructed are particularly interesting
because of the role they play in the structure theory of Delaunay
polytopes and Delaunay tilings, which are described in the following
section. \

\section{Homogeneous and Inhomogeneous \\ Forms on Lattices\label{forms}}

The Voronoi and Delaunay tilings for point lattices are constructed using
the Euclidean metric, but are most effectively studied by mapping the
lattice onto $\mathbb{Z}^{d}$, and replacing the Euclidean metric by an
equivalent metrical form. \ \ For a $d$-dimensional point lattice $\mathbf{%
\Lambda }$ with a basis $\mathbf{b}_{1},\mathbf{b}_{2},...,\mathbf{b}_{d}$
this is done as follows. \ \ If $\mathbf{v}$ is a lattice vector with
coordinates $z_{1},z_{2},...,z_{d}$ relative to this basis, then $\mathbf{%
v}$ can be written as $\mathbf{v}=\mathbf{Bz}$, where $\mathbf{%
B=[\mathbf{b}_{1},\mathbf{b}_{2},...,\mathbf{b}_{d}]}$ is the basis
matrix.
\ \ The squared Euclidean length, is given by $|\mathbf{v }|^{2}=%
\mathbf{z}^{T}\mathbf{B}^{T}\mathbf{Bz=}\varphi
_{\mathbf{B}}(\mathbf{z})$,
where $\mathbf{z}$ is the column vector given by $\mathbf{[}%
z_{1},z_{2},...,z_{d}]^{T}$. \ This squared length can equally well be
interpreted as the squared length of the integer vector $\mathbf{z\in }%
\mathbb{Z}^{d}$ relative to the metrical form $\varphi _{\mathbf{B}}$. \
Therefore, the Voronoi and Delaunay tilings for $\mathbf{\Lambda }$,
constructed using the Euclidean metric, can be studied using the
corresponding Voronoi and Delaunay tilings for $\mathbb{Z}^{d}$
constructed using the metrical form $\varphi _{\mathbf{B}}$. \ \ Moreover,
variation of the Voronoi and Delaunay tilings for $\mathbf{\Lambda }$ in
response to
variation of the lattice basis $\mathbf{b}_{1},\mathbf{b}_{2},...,\mathbf{b}%
_{d}$ can be studied by varying the metrical form $\varphi $ for the fixed
lattice $\mathbb{Z}^{d}$. \ \ In the discussion below we will keep the
lattice fixed at $\mathbb{Z}^{d}$, and vary the metrical form $\varphi $.
With slight abuse of terminology we call any semidefinite form
\emph{metric}.

\begin{definition}
The \emph{squared length} of a vector $\vv$ with respect to a metric form
$\varphi$ is called \emph{the norm of $\vv$ relative to $\varphi$.}  \ \
\end{definition}

\noindent \textbf{The inhomogeneous domain of a Delaunay polytope: } The
following discussion of the geometry of $\cP$ is a summary of some results
contained in (Erdahl, 1992). \ Let $\cP$ be the cone of quadratic
functions on $\mathbb{R}^{d}$ defined by:
\begin{equation*}
\cP = \{\, f \in \R [x_1,\ldots,x_n] \; \vline \; \deg f = 2, \;\; \forall
\mathbf{z} \in \mathbb{Z}^{d} \: f(\mathbf{z}) \geq 0 \, \}.
\end{equation*}
The condition $\forall \mathbf{z}\in \mathbb{Z}^{d} \;\; f(\mathbf{z})\geq
0$ requires  the quadratic part of $f$ to be positive semi-definite, and
requires  any subset of $\mathbb{R}^{d}$ where $f$ assumes negative values
to be free of
$\mathbb{Z}^{d}$-elements. \ The real surface determined by the equation $f(%
\mathbf{x})=0$ might be empty set, it might be a subspace, or it might
have the form
\begin{equation*}
\mathcal{E}_{f}=\mathcal{E}_{0}\times K,
\end{equation*}
where $\mathcal{E}_{0}$ is an ellipsoid and $K$ a complementary subspace.
\
The last case is the interesting one - depending on the dimension of $K$, $%
\mathcal{E}_{f}$ is either an empty (of lattice points) ellipsoid or an
empty cylinder with ellipsoid base. \

We denote by $V(f)$ the set $\mathcal{E}_{f}\, \cap\, \mathbb{Z}^{d}$. In
the case where the surface $\mathcal{E}_{f}$ includes integer points and
is an empty ellipsoid, $V(f)=\mathcal{E}_{f}\, \cap\, \mathbb{Z}^{d}$ is
the vertex set for the corresponding Delaunay polytope $D_{f} = \conv
V(f)$. \
Conversly, if $D$ is a Delaunay polytope in $\mathbb{Z}%
^{d}$, there is a circumscribing empty ellipsoid $\mathcal{E}_{D}$\
determined by a function $f_{D}\in \cP$. \ More precisely, since $D$ is
assumed Delaunay there is a metrical form $\varphi _{D}$, a center $\mathbf{c%
}$ and radius $R$, so that $f_{D}(\mathbf{x})=\varphi_{D}(\mathbf{x}-%
\mathbf{c})-R^{2}=0$ is the equation of a circumscribing empty ellipsoid $%
\mathcal{E}_{D}$. \ Since  $f_{D}$ is non-negative on $\mathbb{Z}^{d}$, it
is an element of $\cP$.

In the case where $\mathcal{E}_{f}$ is an empty cylinder, there will  be
an infinite number of integer points lying on this surface. \ When this
happens $V(f)~=~\mathcal{E}_{f}~\cap~\mathbb{Z}^{d}$ is the vertex set for
a non-bounded \textit{Delaunay polyhedron }$D_{f}=\conv V(f)$, which is a
cylinder with a base, which is also called a Delaunay polytope.  Recall
than a convex polyhedron is called a polytope when it is bounded.\

\begin{definition}
Let $D$ be a Delaunay polyhedron. \ Then, the (inhomogeneous) domain
$\cP_{D}$\ for $D$ is:
\begin{equation*}
\cP_{D}=\{\,f \in \cP \;\; \vline \;\;\text{ }D_{f}=D \,\}
\end{equation*}
\end{definition}

\noindent Such domains are relatively open convex cones that
partition~$\int \cP$ (in fact, they partition a larger subset of $\cP$,
which includes $\int \cP$).

The elements $f\in \cP_{D}$ satisfy the homogeneous linear equations
$f(\mathbf{z})~=~0,\; \mathbf{z} \in \vert~D = D~\cap~\mathbb{Z}^{d}$, so
the dimensions of domains vary depending upon the rank of this system. \ \
\ When $D$ is a single edge, the rank is one and $\cP_{D}$ is a relatively
open facet of the partition with dimension $\dim \cP-1=$
$\binom{d+2}{2}-1$. \ When the rank is equal to $\binom{d+2}{2}-1$, which
is the maximum possible in order that $\cP_{D}$ not be empty, $\cP_{D}$ is
an extreme ray of $\cP$. \

\begin{definition}
Let $V(p)$ be the set of integer points lying on the boundary of a
$d$-dimensional Delaunay polyhedron. A function $p\in \cP$ is perfect if
and only if the equations $p(\mathbf{z})~=~0, \quad \mathbf{z}\in V(p)$
uniquely determine $p$ up to scaling. \ In this case we will call the
subset $V(p) \subset \mathbb{Z}^{d}$ perfect, and will call $D_{p}=\conv
V(p)$ perfect.
\end{definition}

\noindent The elements of a 1-dimensional inhomogeneous domain are
perfect, and the Delaunay polytopes that determine such domains are
perfect. \ The perfect subsets $V(p)$ must be maximal among the subsets
$\{\, V(f) \subset \mathbb{Z}^{d} \;\; \vline \;\; f\in \cP \,\}$. \

\ Perfect inhomogeneous forms are analogues of perfect homogeneous
forms---both achieve their minimum value on $\mathbb{Z}^{d}$ sufficiently
often that the representations of the minimum determine the form. \ In the
homogeneous case the minimum is taken on the non-zero elements of
$\mathbb{Z}^{d}$; if the minimum value $m_{\varphi }$ is achieved on the
set \ $V(
{\varphi })$, then $\varphi $ is uniquely determined by the equations $%
\varphi (\mathbf{z})=m_{\varphi }, \quad \mathbf{z\in }V({\varphi})$. \
Similarly, in the inhomogeneous case the minimum is taken on $\mathbb{Z%
}^{d}$, and the minimum value is zero; if this minimum value is achieved
on the set \ $V(f)\subset \mathbb{Z}^{d}$, then $f$ is uniquely
determined (up to a scale factor) by the equations $f(\mathbf{z})=0, \quad \mathbf{%
z\in }V(f).$ \

\begin{theorem}
If $p\in \cP$ is perfect, then there is an ellipsoid $\mathcal{E}_{0}$,
where $0\leq
\dim (\mathcal{E}_{0})\leq d$, and complementary subspace $K$ so that $%
\mathcal{E}_{p}=\mathcal{E}_{0}\times K$, and so that $\mathcal{E}_{0}$ and $%
K$ satisfy the following arithmetic conditions:
\begin{eqnarray*}
&&\mathcal{E}_{0}\cap \mathbb{Z}^{d}\text{ are the vertices of a perfect
Delaunay polytope in }(\aff~\mathcal{E}_{0}) \cap  \mathbb{Z}^{d}\text{;} \\
&&K\cap \mathbb{Z}^{d}\text{ is a sublattice such that }\mathbb{Z%
}^{d}=((\aff~\mathcal{E}_{0}) \cap  \mathbb{Z}^{d}) \oplus (K\cap
\mathbb{Z}^{d}).
\end{eqnarray*}
In this case we have $V(p)=\mathcal{E}_{p}\cap
\mathbb{Z}^{d}=(\mathcal{E}_{0}\cap \mathbb{Z}^{d})\bigoplus (K\cap \mathbb{Z}%
^{d})$. \ Conversely, any ellipsoid $\mathcal{E}_{0}$ and an affine
complementary
subspace $K$, satisfying these arithmetic conditions, determine a surface $%
\mathcal{E}_{0}\times K$ for a perfect element $p$\ of $\cP$ and,
therefore, determine a perfect element up to a scale factor. \ \
\end{theorem}

\noindent By this theorem, if $p$ is perfect, then $D_{p}$ is either a
perfect Delaunay polytope or, in the degenerate case when $\dim K>0$, a
cylinder with a perfect Delaunay polytope as base. \

A $\Z^d$-vector $[u_1,\ldots,u_d]$ is called \emph{primitive} if
$\GCD(u_1,\ldots,u_d)=1$. As an example, consider primitive vectors $
\mathbf{a,b} \in \mathbb{Z}^{d}$ \
such that $\mathbf{a\cdot b=}1$. \ Then $\mathcal{E}_{0}=\{%
\mathbf{0,b} \}$ is an empty ellipsoid in the 1-dimensional lattice $(\aff%
\mathcal{E}_{0})\cap \mathbb{Z}^{d}$, and $K=\mathbf{a}^{\perp } = \left\{
\x \in \R^d \: | \: \x \cdot \mathbf{a} = 0 \right\}$ is a complementary
subspace. The subset $\{ \mathbf{0,b} \}= \mathcal{E}_{0} \cap
\mathbb{Z}^{d}$ is also the vertex set for a Delaunay polytope in $\aff
(\mathcal{E}_{0} \cap \mathbb{Z}^{d})$, and $K$ is defined so that
$\mathbb{Z}^{d}=((\aff \mathcal{E}_{0}) \cap \mathbb{Z}^{d}))\bigoplus
(K\cap \mathbb{Z}^{d})$. \ These are the arithmetic conditions stated in
the theorem, so that $\{\mathbf{0,b} \} \times K$ is the surface for a
perfect element $p$\ of $\cP$, which, up to a scale factor, is given by
$p(\mathbf{x})=(\mathbf{a\cdot x})( \mathbf{a \cdot x-}1)$. \ The
preimages of negative real values of $p$ lie between two hyperplanes with
equations $\mathbf{a\cdot x=}0$ and $\mathbf{a\cdot x=}1 $, a region which
is a degenerate Delaunay polyhedron.

With the exception of the one-dimensional perfect domains, all
inhomogeneous domains $\cP_{D}$ have proper faces that are inhomogeneous
domains of lesser dimensions. \ If $D$ is inscribed into another
$\Z^d$-polytope $D^{\prime }$ so that $\vert D^{\prime} \supset \vert D$,
then $\cP_{D^{\prime }}\subset
\partial \cP_{D}.$ \

In the following definition the Delaunay polyhedron $D$ may be of any
dimension.

\begin{definition}
Let $D$ be a lattice Delaunay polyhedron. \ Then, the \emph{inhomogeneous
domain} for $D$ is defined as:
\begin{equation*}
\cP_{D}=\{\, f\in \cP \; \; \vline  \;\; \vert D = V(f) \,\}
\end{equation*}
\end{definition}
\noindent If $D$ is a $d$-dimensional Delaunay simplex then $\dim \cP_{D}=$ $%
\binom{d+1}{2}$, but if $D$ is a perfect Delaunay polyhedron, then $\dim
\cP_{D}=1$.

\begin{theorem}
Let $D$ be a $d$-dimensional Delaunay polyhedron.  Then
\begin{equation*}
\cP_{D}=\{\, \sum_{\{p\; \; | \;\; V(p)\supseteq \vert D\}}\omega _{p}p
\;\; \vline \;\; \omega_p \in \R_{>0} \,\},
\end{equation*}
where the summation is over all perfect elements $p$ such that $ \vert D
\subseteq V(p)$.
\end{theorem}

In general, not all relatively open faces on the boundary of an
inhomogeneous domain $\cP_{D}$ are inhomogeneous domains. \ However, in
the special case where $D$ is $d$-dimensional, all extreme rays are
perfect inhomogeneous domains and all relatively open faces are
inhomogeneous domains (see Erdahl 1992). \ In this case an arbitrary
element $f\in \cP_{D}$ has the following representation:
\begin{equation*}
f=\sum_{\{p\;\; \vline \;\; V(p) \supseteq \vert D\}}\omega_{p}p,
\end{equation*}
where $\omega _{p} >0.$ \ \ The summation is over the \emph{perfect}
elements $p$ with the property that $\vert D \subseteq V(p).$  \

These last two paragraphs, and our representation theorem for Delaunay
polytopes, show the important role played by perfect Delaunay polytopes in
the theory. \ \ \ \

\medskip

\noindent \textbf{The homogeneous domain of a Delaunay tiling: }\
Voronoi's classification theory for lattices, his \textit{theory of
lattice types} (L-types), was formulated using metrical forms and the
fixed lattice $\mathbb{Z}^{d}$. \ In this theory two lattices are
considered to be the same type if their Delaunay tilings are affinely
equivalent. \ Consider a positive definite quadratic form $\varphi $. \
Then a lattice polytope $D$ is Delaunay relative to $\varphi $ if it can
be circumscribed by so called \emph{empty ellipsoid }$\mathcal{E}$ with
equation of the form
\begin{equation*}
\varphi (\mathbf{x-c})=R^{2}\text{,}
\end{equation*}
where  $\mathbf{c}\in \mathbb{R}^{d}$ and  $R \in \R_{>0}$; in addition,
 $\conv \mathcal{E}$ must have no interior $\mathbb{Z}^{d}$%
-elements, and $\vert D$ must be given by $\mathcal{E\cap }%
\mathbb{Z}^{d}$. \ The collection of all such Delaunay polytopes fit
together facet-to-facet to tile $\mathbb{R}^{d}$, a tiling that is
uniquely
determined by $\varphi $. \ This is the Delaunay tiling $\mathcal{D}%
_{\varphi }$ for $\mathbb{Z}^{d}$ relative to the metrical form $\varphi
$.
\ If a second metrical form $\vartheta $ has Delaunay tiling $\mathcal{D}%
_{\vartheta }$, and if $\mathcal{D}_{\vartheta }$ is identicle to, or
affinely equivalent to $\mathcal{D}_{\varphi }$, then $\varphi $ and $%
\vartheta $ are metrical forms of the same L-type. \

The description we give below requires that certain degenerate metrical
forms be admitted into the discussion, namely, those forms $\varphi $ for
which $K=\Ker \varphi$ is a rational subspace of $\R^d.$ \ \ \ The
Delaunay polyhedra for such a form are themselves degenerate--they are
cylinders with axis $K$ and Delaunay polytopes
as bases. \ These cylinders fit together to form the degenerate Delaunay tiling $%
\mathcal{D}_{\varphi }$. \ \ For example, if $\mathbf{a\in
}\mathbb{Z}^{d}$ is primitive, $\varphi (\mathbf{x})=(\mathbf{a\cdot
x})^{2}$ is such a form---the kernel $K$ is the solution set for $\varphi
(\mathbf{x)=}0$, and given by $\mathbf{a}^{\perp }$, which is rational. \
The Delaunay tiles are
infinite slabs, each bounded by a pair of hyperplanes $\mathbf{a\cdot x=}k%
\mathbf{,a\cdot x=}k+1$, $k\in \mathbb{Z}$. \ These fit together to tile $%
\mathbb{R}^{d}$. \

Let $\Phi $ be the cone of metrical forms in $d$ variables, namely, the
cone of positive definite quadratic forms and semidefinite quadratic forms
with rational kernels. \ For
each metrical form $\varphi \in \Phi $ there is a Delaunay tiling $\mathcal{D%
}_{\varphi }$ for $\mathbb{Z}^{d}.$

\begin{definition} If $\mathcal{D}$ is a Delaunay tiling for $\mathbb{Z}^{d}$, then
the following cone of positive definite quadratic forms
\begin{equation*}
\Phi _{\mathcal{D}}=\left\{\, \varphi \in \Phi^d \;\; \vline  \;\; \mathcal{D}_{\varphi }=%
\mathcal{D}\, \right\}
\end{equation*}
is called an L-type domain.
\end{definition}

\noindent For this definition the Delaunay tilings can be the usual ones,
where the tiles are Delaunay polytopes -- or they could be degenerate
Delaunay tilings where the tiles are cylinders with a common axis $K$.

The relatively open faces of an L-type domain, are L-type domains. \ If
$\mathcal{D}$ is a trangulation, $\Phi _{ \mathcal{D}}$ has full dimension
$\binom{d+1}{2}$; this is the generic case. \ If $\mathcal{D}$ is not a
triangulation, $\dim \Phi _{\mathcal{D}}$ is less than $\binom{ d+1}{2}$,
and $\Phi _{\mathcal{D}}$ is a boundary cell of a full-dimensional L-type
domain. \ In the case where $ \dim \Phi _{\mathcal{D}}=1$, the elements
$\varphi \in \Phi _{\mathcal{D}}$ are called edge forms as they correspond
to extreme rays of full-dimensional L-type domains.

\par
Let $D$ be a  polytope in the Delaunay
tiling $\mathcal{D}$. \ Then, if $\pi _{\Phi}$ is the projection onto the quadratic part, $%
\Phi _{\mathcal{D}}\subset \pi _{\Phi}(\cP_{D})$. \ \ Since this
containment holds for all Delaunay tiles $D\in \mathcal{D}$, there is the
following description of $\Phi _{\mathcal{D}}$ in terms of inhomogeneous
domains:
\begin{equation*}
\Phi _{\mathcal{D}}=\bigcap_{D\in \mathcal{D}}\pi _{\Phi}(\cP_{D}),
\end{equation*}
where the intersection is over all $d$-dimensional Delaunay polytopes in $
\mathcal{D}$. \ Since the containment $\Phi_{\mathcal{D}}\subset
\pi_{\Phi}(\cP_{D})$ also holds for all Delaunay tilings $\mathcal{D}$
that contain $ D$, there is the following description of
$\pi_{\Phi}(\cP_{D})$ in terms of homogeneous domains:
\begin{equation*}
\pi_{\Phi}(\cP_{D})=\bigsqcup_{\mathcal{D\ni }D}\Phi _{\mathcal{D}},
\end{equation*}
where the disjoint union is over all Delaunay tilings for $\Z^d$ that
contain $D$. \ This holds not only for full-dimensional  cells of
$\cP_{D}$, but for cells of all dimensions. \ The last equality shows that
$\pi _{\Phi}(\cP_{D})$ is tiled by L-type domains. It also establishes the
following

\begin{theorem}
If $p\in \cP$ is perfect and, therefore, an edge form for an inhomogeneous
domain, then $\pi _{\Phi }(p)$ is an edge form for an L-type domain. \ \
\end{theorem}

\noindent This result can also be established in a more direct way by
appealing to the definition of L-type domain: if $D$ is a perfect Delaunay
polytope, or perfect degenerate Delaunay polyhedron, then $\pi _{\Phi
}(\cP_{D})$ is a one-dimensional L-type domain. \ By definition, the
elements of $\pi _{\Phi }(\cP_{D})$ are then edge forms. \

The converse of theorem does not hold - there are edge forms for L-type
domains that are not inherited from perfect elements in $\cP$. \ \
Evidence has accumulated recently that the growth of numbers of types of
edge forms with dimension is very rapid starting in six dimensions (see
Dutour and Vallentin, 2003), but the the growth of perfect inhomogeneous
forms is much less rapid - there is some hope that a complete
classification can be made through dimension nine, or even through ten. \

\section{Infinite Sequences of \\ Perfect Delaunay polytopes}

Consider the following sets of vectors in $\mathbb{R}^{d}$:
\begin{equation*}
D_{s,k}^{d}=\{\, [1^{s},0^{d-s}]-\frac{s-1}{d-2k} \mathbf{j}\,\}\cup
\{\,[1^{s+1},0^{d-s-1}]-\frac{s}{d-2k}\mathbf{j}\,\}
\end{equation*}
where $s,k \in \N$ and $\mathbf{j=[}1^{d}]$ ($1^{d}$ means
the entry $1$ is repeated $d$ times, and similarly for $1^{s}$ and $0^{d-s}$%
). \ All permutations of entries are taken so that $|D_{s,k}^{d}|=\binom{d}{s%
}+\binom{d}{s+1}=\binom{d+1}{s}$. The following is the main theorem of
this paper.\ \

\begin{theorem}\label{main_theorem}
For $d\geq k(2s+1)+1$, where $s\geq 1,\;k\geq 2$, the polytope
\begin{equation*}
P_{s,k}^{d}=\frac{1}{2}\conv\{\;D_{s,k}^{d}\cup -D_{s,k}^{d}\;\}.
\end{equation*}
is a symmetric perfect Delaunay polytope for the affine lattice $\Lambda _{s,k}^{d}=\aff C_{s,k}^{d}$; the origin is the center of symmetry for $%
P_{s,k}^{d}$ and does not belong to $\Lambda _{s,k}^{d}$. \ The
circumscribing empty ellipsoid can be defined by the equation $\varphi _{s,k}^{d}(%
\mathbf{x})~=~R^{2}$, where
\begin{equation*}\label{formula_for_phi}
\varphi
_{s,k}^{d}(\mathbf{x})=4k(d-k(2s+1))|\mathbf{x}|^{2}+(d^{2}-(4k+2s+1)d+4k(2s+k))(\mathbf{j\cdot
x)}^{2}.
\end{equation*}
\end{theorem}

\noindent Each pair of positive \ integers $(s,k)$, for $s\geq 1,k\geq 2$,
determines an infinite sequence of symmetric perfect Delaunay polytopes,
one in each dimension, with the initial dimension given by $k(2s+1)+1$.
For $s=1,k=2$ the infinite sequence is the one described in the opening
commentary, i.e.,   C$^{d},d\geq 7$, where the initial term is the Gosset
polytope $3_{21}=\textrm{C}^7$ with $56$ vertices. \

The $\binom{d+1}{s}$ diagonal vectors $D_{s,k}^{d}$ for $P_{s,k}^{d}$ have
the origin as a common mid-point, forming a segment arrangement that
generalizes the cross formed by the diagonals of a cross-polytope. \
Moreover, these \ $\binom{d+1}{s}$ diagonals are primitive and belong to
the same \textit{parity class} for $\Lambda_{s,k}^{d}$, namely, they are
equivalent modulo $2\Lambda_{s,k}^{d}$. \
 \ \ More generally, primitive
lattice vectors $ \uu ,\vv $ in some lattice $\Lambda $, with mid-points
equivalent modulo $\Lambda$, are necessarily equivalent modulo $2\Lambda $
and thus belong to the same parity class. \ \ And conversely, the
mid-points of lattice vectors $\uu ,\vv \in \Lambda $ belonging to the
same parity class are equivalent modulo $\Lambda $. By analogy with the
case of cross-polytopes, we call any such arrangement of segments or
vectors a cross. The convex hulls of such crosses often appear as cells in
Delaunay tilings -- cross polytopes are examples, as are the more
spectacular symmetric perfect Delaunay polytopes.  \ There is a criterion
due to Voronoi (1908) and Baranovskii (1991) that determines whether these
crosses are Delaunay:  \emph{Let $\Lambda $ be a lattice, let $\varphi $
be a metrical form, and let $C$  be the convex hull of a cross of
primitive vectors belonging to the same parity class. \ Then $C$ is
Delaunay relative to $\varphi $ if and only if the set of diagonal vectors
forming the cross is the complete set of vectors of minimal length,
relative to $\varphi$, for their parity class. } \ \ This is the criterion
we have used to establish the Delaunay property for the symmetric perfect
Delaunay polytopes $P_{s,k}^{d}.$

The following result shows that asymmetric perfect Delaunay polytopes can
appear as sections of symmetric ones.

\begin{theorem}
For $d\geq 6$ let $\mathbf{u}=[-1^{2},1^{d-1}] \in \Z^{d+1}$. \ Then
\begin{equation*}
\textrm{G}^{d}=\conv\{\, \mathbf{v} \in \vert P_{1,2}^{d+1} \;\; |  \;\;
\mathbf{v \cdot u}=\frac{1}{2} \, \}
\end{equation*}
is an asymmetric perfect Delaunay polytope for the affine sublattice
$\Lambda _{s,k}^{d}= \aff \vert \textrm{G}^{d}$. \ \ The circumscribing
empty ellipsoid can be defined as $\{  \, \x \in \aff \textrm{G}^{d} \;\;
| \;\; \varphi _{1,2}^{d+1}(\mathbf{x})=R^{2} \, \}$, where
\begin{equation*}
\varphi _{1,2}^{d+1}(\mathbf{x})=8(d-5)|\mathbf{x}|^{2}+(d^{2}-9d+22)(\mathbf{%
j\cdot x)}^{2}.
\end{equation*}
\end{theorem}

\noindent This is the infinite sequence of G-topes described in the
introduction, with  Gossett polytope $2_{21}=\textrm{G}^6$ as the initial
term for $d=6$. \

The terms in this sequence have similar combinatorial properties. \ For
example, the lattice vectors running between vertices all lie on the
boundary, in all cases. These lattice vectors are either edges of
simplicial facets, or diagonals of cross polytope facets--there are two
types of facets, simplexes and cross polytopes. \ The 6-dimensional
Gossett polytope has $27$ five-dimensional cross-polytopal facets, but for
the rest the number is twice the dimension, $2d.$ G$^6$ can be found as a
section of G$^7$, but G$^8$ and G$^9$ do not have sections arithmetically
equivalent to G$^6.$

From this theorem we might expect that all asymmetric perfect Delaunay
polytopes are obtained from symmetric ones by sectioning, and at our
current state of knowledge this appears to be the case. \ This would be a
fortunate turn of events, since the growth of classes of symmeteric
perfect Delaunay polytopes with dimension is slower than for all other
geometric phenomena associated with point lattices that we know about. \ \

We summarize properties of the constructed polytopes in the following
table (in this notation groups of entries separated by semicolumns can
only be permuted between themselves).

\begin{table}
\[ \begin{array}{|c|c|c|c|} \hline
P & \dim P & | \vert P| & Symmetry \\ \hline P_{\kappa,s}^d & d &
2\left({d~\choose~s}+{d\choose{s+1}}\right)=2{{d+1}\choose {s+1}} &
\textrm{centrally-symmetric} \\ \hline \Upsilon^{d-1} & d-1 &
\frac{d(d+1)}{2} - 1 & asymmetric \\ \hline
\end{array} \]

\caption{Properties of constructed perfect Delaunay polytopes}
\label{degree-4-number-table}
\end{table}

The following table gives coordinates of the vertices of G-topes,
discovered in 2001 by Erdahl and Rybnikov.
\begin{center}
\begin{tabular}{|c|c|c|}
\hline
$\lbrack 0^{d}]\times 1$ & $[-1,0^{d-1};-1]\times (d-1)$ & $%
[1^{d-1};-(d-3)]\times (d-1)$ \\ \hline
$\lbrack 0,1^{d-2};(d-4)]\times (d-1)$ & $[1^{2},0^{d-3};1]\times \frac{%
(d-1)(d-2)}{2}$ & $[1,0^{d-1};0]\times (d-1)$ \\ \hline
\end{tabular}
\end{center}

Polytope $P_{4,1}^7$ is affinely equivalent to Gosset's $3_{21}=\textrm{
G}^7$. Polytope $\Upsilon^6$ is affinely equivalent to Gosset's
$2_{21}=\textrm{ G}^6$.

\bigskip

\section{Delaunay Property Part of Main Theorem} We will use the
following notation: $\phi_1 = \left(\sum_{i=1}^d x_i\right)^2$,
$\phi_2(\x) = \left| \x - {\bf j}\frac{\sum_{i=1}^d x_i}{d} \right|^2$;
$\phi_2(\x)$ is the squared Euclidean distance from $\x$ to the line $\lin
{\bf j}$. The following theorem proves the Delaunay property of polytopes
$P_{s,k}^d$ asserted in Theorem 4.

\begin{theorem}
\label{delaunay-property-theorem} Let $s\geq 1,k\geq 2$, and  $d\geq
k(2s+1)+1$. Let $l^0 = [{(-1)}^{k/2}, 1^{d-{k/2}}]$, $\Lambda =
{<e_1,\dots,e_d,\frac{\bf j}{d-k}>}_{\mathbb Z}$, $\Lambda^0 = \lbrace
\lambda \in \Lambda: \lambda \cdot l^0 \equiv 1 (\pmod 2) \rbrace$. Then
there is a positive definite quadratic form of the type
\begin{equation}\label{quadratic_form_equation} \phi(\x) = \alpha \phi_1(\x) + \beta
\phi_2(\x),
\end{equation}
were $\alpha, \beta >0$, such that the polytope $P_{s,k}^d$ is a Delaunay
polytope in the affine lattice $\Lambda^0$ with respect to $\phi$.
\end{theorem}

Here $n=d - k$ and $e_1,\dots,e_d$ stand for the canonical basis of
${\mathbb R}^d$.

\begin{lemma} \label{structural_lemma}
Suppose $\phi(\x) = \alpha \phi_1(\x) + \beta \phi_2(\x)$ where $\alpha,
\beta > 0$. Then all points $\lambda\in\Lambda^0$ which are closest to $0$
with respect to $\phi$, ie.,
$$ \phi(\lambda) = \min\lbrace \phi(u) : u\in\Lambda^0 \rbrace $$
are, up to permutations of components, of the type $[ 1^l, 0^{d-|l|} ] +
a\frac{{\bf j}}{n}$, where $-\frac{d}{2} \le l < \frac{d}{2}$,
$a\in{\mathbb Z}$. Here $1^l$ for a negative $l$ means $l$ times $-1$.
This representation is unique.
\end{lemma}

\begin{proof} Suppose $\lambda\in \Lambda^0$ is closest to $0$ with respect to
$\phi$ in $\Lambda^0$. Let $\lambda = a_1 e_1 + \dots + a_d e_d + a
\frac{\bf j}{n}$, $a_1,\dots,a_d,a \in {\mathbb Z}$, $m = a_1 + \dots +
a_d$. We have
\begin{equation}
\phi_2(\lambda) = \phi_2\left( \sum_{i = 1}^{d} a_i e_i \right) =
 \left| \sum_{i = 1}^{d} a_i e_i - {\bf j}\frac{m}{d} \right|^2 = \sum_{i=1}^{d} a_i^2 - \frac{m^2}{d}.
\end{equation}

The numbers $a_i$ are of two consecutive integer values. Otherwise there
would be $a_i$ and $a_j$ such that $a_i - a_j \ge 2$. Consider the vector
$\lambda' = a_1 e_1 + \dots + (a_i - 1) e_i + \dots + (a_j + 1) e_j +
\dots + a_d e_d + a \frac{\bf j}{n}$. We have $\phi_2(\lambda') -
\phi_2(\lambda) = (a_i - 1)^2 + (a_j + 1)^2 - a_i^2 - a_j^2 = 2(a_j - a_i)
+ 2 \le -2$ and ${\bf j}\cdot \lambda' = {\bf j}\cdot \lambda$. Since
$\lambda' \in \Lambda^0$ and $\phi(\lambda') < \phi(\lambda)$, it follows
that the vector $\lambda$ is not closest to $0$ which is a contradiction.

Now let $b$ be the smallest of the values of numbers $a_i$. Subtract
$b(e_1 + \dots + e_d)$ from the first part and add an equal value of
$bk\frac{\bf j}{n}$ to the second part of the existing representation of
$\lambda$. After a permutation of components, the vector $\lambda$ is
equal to $[1^l,0^{d-l}] + (a + bk)\frac{\bf j}{n}$ where $0\le l < d$. If
$l \ge \frac{d}{2}$, again subtract $[1^d]$ from the first summand and add
${\bf j}$ to the second summand to get the required representation.

To prove uniqueness, note that the components $\lambda_i$ of the vector $[
1^l, 0^{d-|l|} ] + a\frac{{\bf j}}{n}$ are of at most two values. In
vector $[ 1^l, 0^{d-|l|} ]$, either $1$s fill the positions of the  bigger
$\lambda_i$, or $-1$s fill the positions of the smaller $\lambda_i$. Since
$-\frac{d}{2} \le l < \frac{d}{2}$, the choice is unique. $\Box$
\end{proof}

\begin{lemma}
\label{2nd_structural_lemma} All vectors in the affine lattice $\Lambda^0$
which are closest to $0$ with respect to a form $\phi = \alpha \phi_1(\x)
+ \beta \phi_2(\x)$, $\alpha,\beta > 0$, belong to the set $\lbrace
\lambda \in \Lambda^0 : \left| \lambda \cdot {\bf j} \right| \le
\frac{d}{n}\rbrace$. If $\lambda$ is closest to $0$ with respect to $\phi$
and  $\left| \lambda \cdot {\bf j} \right| = \frac{d}{n}$, then
$\lambda\in\lbrace\pm\frac{\bf j}{n}\rbrace$.
\end{lemma}

\begin{proof} By multiplying both $\alpha$ and $\beta$ by the same positive number, we can assume that
the minimal value of $\phi$ on $\Lambda^0$ is $1$. Since $\frac{\bf j}{n}
\in \Lambda^0$, $\phi(\frac{\bf j}{n}) = \alpha (\frac{d}{n})^2 \ge 1$.
Let $\lambda\in{\mathbb R}^d$ be a point with $\phi(\lambda)=1$. Represent
$\lambda = \frac{\gamma}{n} {\bf j} + u$ where $\gamma\in{\mathbb R}$,
$u\cdot{\bf j} = 0$. We have $\phi(\lambda) = \alpha\frac{\gamma^2
d^2}{n^2} + \beta |u|^2 = 1$, therefore $\alpha\frac{\gamma^2 d^2}{n^2}
\le 1$. Since $\alpha(\frac{d}{n})^2 \ge 1$, we have proven that $|\gamma|
\le 1$, so $|\lambda\cdot{\bf j}| = |\gamma|\frac{d}{n} \le \frac{d}{n}$.

If the latter inequality holds strictly, then $|\gamma|=1$ and
$\phi(\lambda) = \alpha\frac{d^2}{n^2} + \beta |u|^2 = 1$. Since $\alpha
(\frac{d}{n})^2 \ge 1$, we necessarily have $\beta = 0$ so $\lambda =
\pm\frac{\bf j}{n}$.
\end{proof}

\medskip

\par \noindent \textbf{Proof of Theorem 6:}

We define a collection of points in $\Lambda^0$ which contains, up to sign
and permutation of components, all vectors of the affine lattice
$\Lambda^0$ which have the minimal distance to $0$ with respect to any
quadratic form $\phi$ of the type (\ref{quadratic_form_equation}). Then we
pick the form $\phi$ so that the vertices of $P^d_{k,s}$ are minimal
vectors and other points of $M$ are not. Below is an implementation of
this program.

\medskip

(A) Consider the set
\begin{equation}
M = \left\lbrace \lambda = [ 1^l, 0^{d-|l|} ] + a\frac{\bf j}{n}
\in\Lambda^0 : 0 \le \lambda\cdot{\bf j} < \frac{d}{n}, -\frac{d}{2} \le l
< \frac{d}{2} \right\rbrace \cup \lbrace \frac{\bf j}{n} \rbrace.
\end{equation}
By the two previous lemmas, for every minimal vector $\lambda\in\Lambda^0$
of the quadratic form $\phi$, either $\lambda$ or $-\lambda$, after
permutation the components if necessary, belongs to this set.

For $\lambda = [ 1^l, 0^{d-l} ] + a\frac{\bf j}{n}$, we have
\begin{equation}
\lambda\cdot{\bf j} = l + \frac{ad}{n}, \\
\phi_1(\lambda) = \left(l + \frac{ad}{n}\right)^2, \\
\phi_2(\lambda) = |l| - \frac{l^2}{d},\\
l^0\cdot\lambda\equiv a+l \pmod 2.
\end{equation}
From the calculations above, we get
\begin{equation}
\label{M_equation} M = \left\lbrace [ 1^l, 0^{d-|l|} ] + a\frac{\bf j}{n}
: l + a \equiv 1 (2), 0 \le lk + ad < d, -\frac{d}{2} \le l < \frac{d}{2}
\right\rbrace \cup \lbrace \frac{\bf j}{n} \rbrace.
\end{equation}

We define a mapping from ${\mathbb R}^d$ to ${\mathbb R}^2$ by
$\phi_{1,2}(\x)=(\phi_1(\x),\phi_2(\x))$. We will call the image of $M$
the {\em diagram}. Lines parallel to $\phi_1$-axis in ${\mathbb R}^2$ will
be called horizontal.

\medskip

(B) For $\lambda_1,\lambda_2\in M$, points
$\phi_{1,2}(\lambda_1),\phi_{1,2}(\lambda_2)$ belong to the same
horizontal line if and only if $\lambda_1 = \pm\lambda_2$. If $\lambda$
and $-\lambda$  both belong to $M$, then $\lambda\cdot {\bf j} = 0$.

To prove that, take $\lambda= [ 1^l, 0^{d-|l|} ] + a\frac{\bf j}{n} \in
M$. If $l\ne 0$, then the condition $0 \le lk + ad < d$  uniquely defines
$a$ from a given value of $l$, and if $l=0$, then $a=1$. Therefore $l$
uniquely defines $a$.

The function $l\to |l| - \frac{l^2}{d}$ is even and increasing on
$[0,\frac{d}{2}]$ so two different points of the diagram may belong to the
same horizontal line if and only if their preimages are $\lambda_1=[ 1^l,
0^{d-|l|} ] + a_1\frac{\bf j}{n}$ and $\lambda_2=[ 1^{-l}, 0^{d-|l|} ] +
a_2\frac{\bf j}{n}$ for some $l$,$a_1$,$a_2$. If $l\ne 0$, then $0 \le lk
+ a_1d < d$, $0 \le -lk + a_2d < d$. Adding these inequalities, we get $0
\le (a_1 + a_2)d < 2d$, so $a_1 + a_2 \in\lbrace 0,1 \rbrace$. If
$a_1+a_2=0$, then $\lambda_1=-\lambda_2$ so $\phi_{1,2}(\lambda_1) =
\phi_{1,2}(\lambda_2)$. If $a_1+a_2=1$, then the numbers $l+a_1$ and $-l +
a_2$ have different parity which contradicts conditions $l+a_1 \equiv -l +
a_2 \equiv 1 (2)$. If $l = 0$, then $a_1 = a_2 = 1$. This proves our
claim.

\medskip

(C) Our method of constructing perfect Delaunay polytopes can be
summarised as follows. We note that each line $\alpha x_1 + \beta x_2 =
1$, where $\alpha,\beta > 0$, which supports the edge of the convex hull
of the diagram, gives rise to a quadratic form $\phi(\x)=\alpha \phi_1(\x)
+ \beta \phi_2(\x)$ such that the ellipsoid $\phi(\x)=1$ circumscribes the
points of the affine lattice $\Lambda^0$ which fall into the endpoints of
the edge and has no lattice points inside.

\medskip

\begin{figure}
\begin{center}
\resizebox{!}{240pt}{\includegraphics[clip=true,keepaspectratio=true]{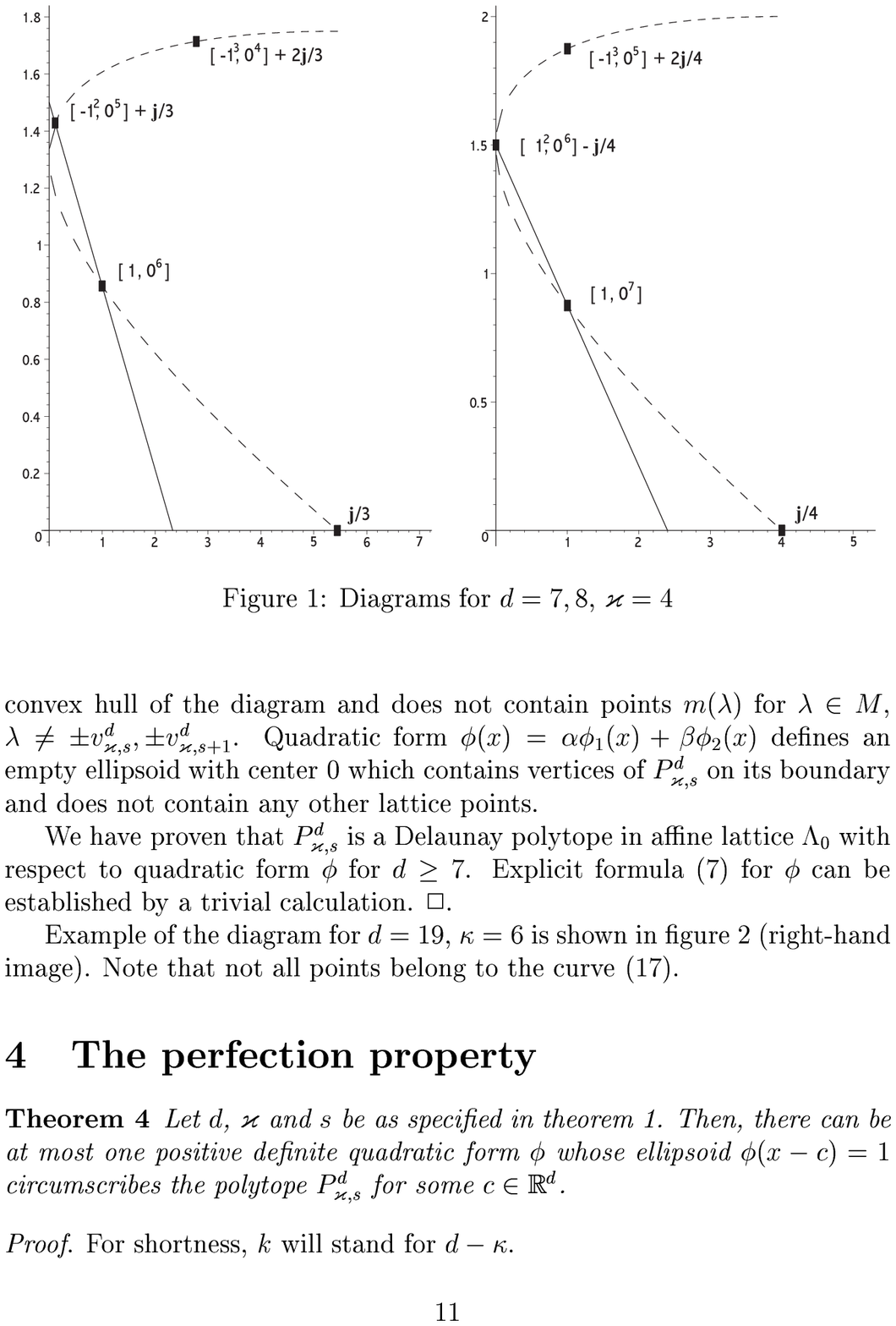}}
\label{7-3-figure}
\end{center}
\end{figure}

(D) First consider case $d=7$, $k=4$, $s = 1$. The diagram for these
values of parameters is shown in figure \ref{7-3-figure}. We see that the
line $\frac{3}{7}x_1 + \frac{2}{3}x_2 = 1$ passes through points
$\phi_{1,2}(v^7_{4,1})=\phi_{1,2}([1,0^6])$ and
$\phi_{1,2}(v^7_{4,2})=\phi_{1,2}([{-1}^2,0^5] + \frac{\bf j}{3})$, and
all other points of the diagram are contained in the open half-plane
$\frac{3}{7}x_1 + \frac{2}{3}x_2 > 1$. This means that polytope
$P_{4,1}^7$ is a Delaunay polytope with respect to quadratic form
$\frac{3}{7}\phi_1(\x) + \frac{2}{3}\phi_2(\x) = \frac{2}{3}\x\cdot \x +
\frac{1}{3}({\bf j}\cdot \x)^2$.

\medskip

(E) Now suppose that $d\ge 8$, $4 \le k \le d/2$, $k\equiv 0(2)$ and $1
\le s \le \frac{d}{k} - 1$. Suppose $\lambda = [ 1^l, 0^{d-l} ] +
a\frac{\bf j}{n}\in M$ and $0 \le \phi_2(\lambda)\le \frac{d}{k} -
\frac{d}{k^2}$ (or, equivalently, $|l| \le \frac{d}{k}$). If $|l| <
\frac{d}{k}$, then $l \ge 0$ and
\begin{equation}
\lambda = [ 1^l, 0^{d-l} ] - (l-1)\frac{\bf j}{n}.
\end{equation}
If $|l| = \frac{d}{k}$ (this implies that $\frac{d}{k}$ is an integer
number), then
\begin{equation}
\pm\lambda = [ 1^\frac{d}{k}, 0^{d-\frac{d}{k}} ] -
(\frac{d}{k}-1)\frac{\bf j}{n}.
\end{equation}
To prove this, first check that these points belong to set $M$, and then
use statement in paragraph (B) that there is at most two points $\lambda=[
1^l, 0^{d-|l|} ] + a\frac{\bf j}{n} \in M$ with a given value of $|l|$,
and there are two if and only if $lk + ad = 0$, in which case the points
are $\pm\lambda$.

We'll use notation $v^d_{k,s} = [ 1^l, 0^{d-l} ] - (l-1)\frac{\bf j}{n}$
for $0 \le s \le \frac{d}{k}$. We have
\begin{equation}
\phi_1(v^d_{k,s}) = {\left(s + (1-s)\frac{d}{n}\right)}^2, \\
\phi_2(v^d_{k,s}) = s - \frac{s^2}{d}. \\
\end{equation}
All points $\phi_{1,2}(v^d_{k,s})$ therefore belong to the curve
\begin{equation}
\label{curve-equation} t \to \left( {\left(t + (1-t)\frac{d}{n}\right)}^2,
t - \frac{t^2}{d} \right)
\end{equation}
which touches the vertical axis in the point $(0,\frac{d}{k} -
\frac{d}{k^2})$ when $t = \frac{d}{k}$. The curve is drawn in dashed line
in figure \ref{7-3-figure}. The portion of the curve for $0 \le t \le
\frac{d}{k}$ is the graph of a convex function.

Therefore, for each $0 \le s \le \frac{d}{k} - 1$ we can find a line with
equation $\alpha x_1 + \beta x_2 = 1$ which passes through points
$\phi_{1,2}(v^d_{k,s})$, $\phi_{1,2}(v^d_{k,s+1})$, supports the convex
hull of the diagram and does not contain points $\phi_{1,2}(\lambda)$ for
$\lambda\in M$, $\lambda \ne \pm v^d_{k,s}, \pm v^d_{k,s+1}$. Quadratic
form $\phi(\x) = \alpha \phi_1(\x) + \beta \phi_2(\x)$ defines an empty
ellipsoid with center $0$ which contains vertices of $P^d_{k,s}$ on its
boundary and does not contain any other lattice points.

We have proven that $P^d_{k,s}$ is a Delaunay polytope in affine lattice
$\Lambda_0$ with respect to a unique quadratic form $\phi=\phi^d_{k,s}$
for $d\ge 7$. Explicit formula (\ref{formula_for_phi}) for $\phi^d_{k,s}$
is  established by a trivial calculation. \begin{flushright} $\Box
$\end{flushright}

\begin{figure}
\begin{center}
\resizebox{!}{240pt}{\includegraphics[clip=true,keepaspectratio=true]{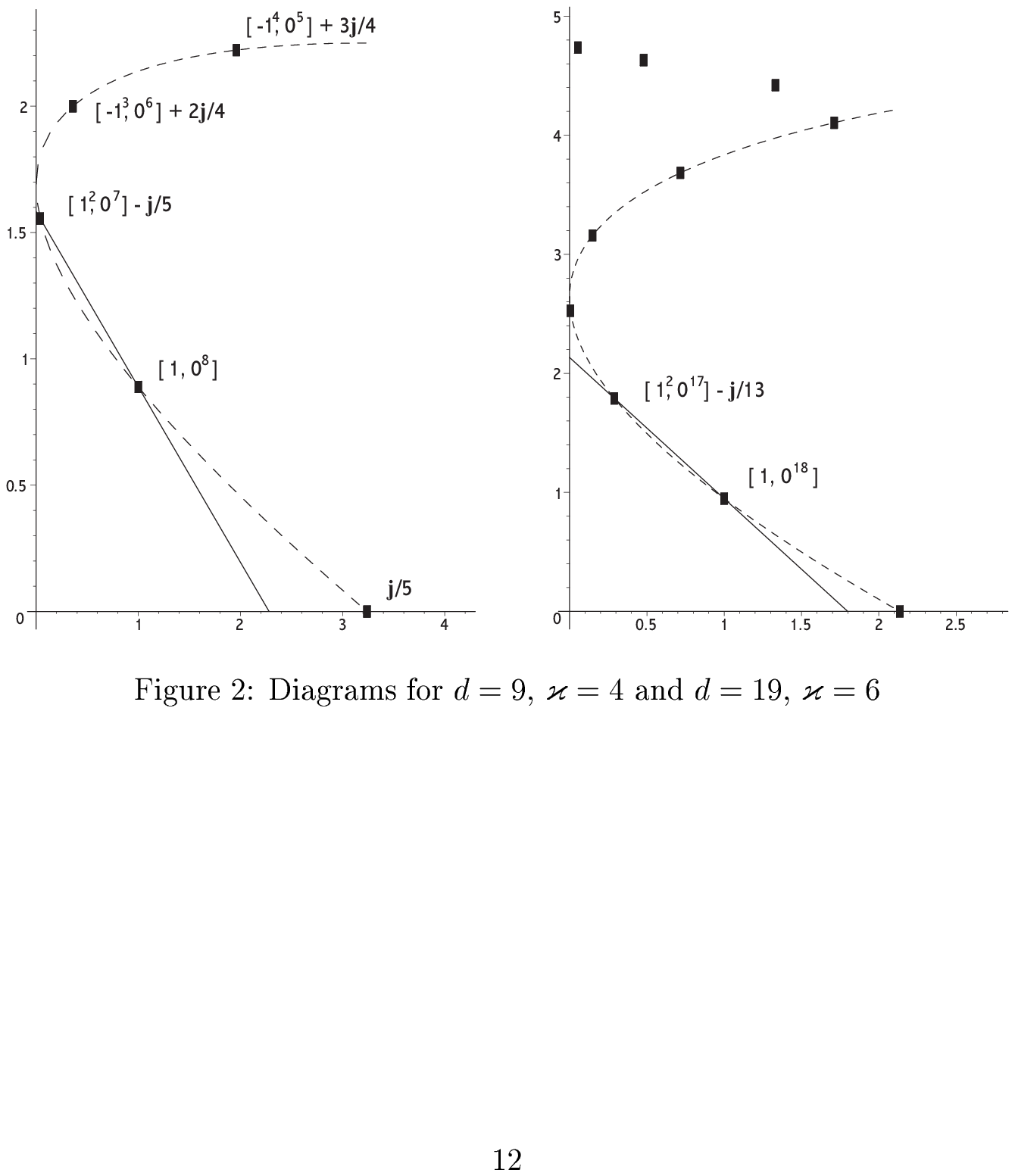}}
\label{9-4-figure}
\end{center}
\end{figure}

Example of the diagram for $d = 19$, $k=6$ is shown in figure
\ref{9-4-figure} (right-hand image). Note that not all points belong to
the curve (\ref{curve-equation}).

\section{Perfection Property Part of Main Theorem}
\begin{theorem}
\label{perfectness-property-theorem} Let $d$, $k$ and $s$ be as in Theorem
\ref{main_theorem}. Then, there can be at most one positive definite
quadratic form $\phi$ whose ellipsoid $\phi(\x-\cc) = 1$ circumscribes the
polytope $P^d_{k,s}$ for some $\cc \in {\mathbb R}^d$.
\end{theorem}

\begin{proof}
As before, $n=d-k$. Since for any given positive definite quadratic form
the center of the its ellipsoid which circumscribes $P^d_{k,s}$ is defined
uniquely, and since the polytope $P^d_{k,s}$ has $0$ as the center of
symmetry, $0$ is the center of the circumscribing ellipsoid. Therefore, if
$\phi$ and $\psi$ are two positive definite quadratic forms which
circumscribe $P^d_{k,s}$ by ellipsoids of radius $1$, then the ellipsoids
are given by equations $\phi(\x) = 1$, $\psi(\x) = 1$. Hence $(\phi -
\psi)|_{V^d_{k,s} \cup V^d_{k,s+1}} = 0$. We consider an arbitrary
quadratic form $f$ such that $f|_{V^d_{k,s} \cup V^d_{k,s+1}} = 0$ and
prove that $f = 0$.

We will use the following symmetrization technique. Suppose that $G$ is a
subgroup of the group $S_d$ of permutations on $d$ elements $I_d = \lbrace
1,\dots, d \rbrace$. The symmetrization $F$ of the form $f$ by $G$ is
defined as
\begin{equation}
F(x_1,\dots,x_d) = \sum_{\sigma\in G} f(x_{\sigma 1},\dots,x_{\sigma d}).
\end{equation}
Since polytope $P^d_{k,s}$ is invariant under transformations of ${\mathbb
R}^d$ defined by permutations of coordinates, all components in the above
sum are $0$ on the set of vertices of $P^d_{k,s}$. Therefore
$F|_{V^d_{k,s} \cup V^d_{k,s+1}} = 0$.

Firstly we prove that if a form $f$ satisfies the condition $f|_{V^d_{k,s}
\cup V^d_{k,s+1}} = 0$, then the sum of the diagonal elements of $f$, and
the sum of off-diagonal elements are 0. We symmetrize $f$ by the group
$S_d$ and get a form $F$ with diagonal coefficients proportional to the
sum of diagonal coefficients of $f$, and non-diagonal coefficients
proportional to the sum of non-diagonal coefficients of $f$. Therefore $F$
can be written as $F(\x) = \alpha \phi_1(\x) + \beta \phi_2(\x)$. Suppose
that $\alpha$ and $\beta$ are not both equal to $0$. Since $F|_{V^d_{k,s}
\cup V^d_{k,s+1}} = 0$, the line $\alpha x_1 + \beta x_2 = 0$ passes
through the points $\phi_{1,2}(V^d_{k,s})$ and $\phi_{1,2}(V^d_{k,s+1})$,
where $\phi_{1,2}(\x)=(\phi_1(\x),\phi_2(\x))$. This is impossible because
the line that passes through these points is uniquely defined and does not
contain $0$. This contradiction proves our claim.

Next we prove that all diagonal elements of the form $f$ are equal to $0$.
Suppose that one of them is nonzero. Without a limitation of generality we
may assume that $t=f_{11} \ne 0$. Let $F$ be the symmetrization of $f$ by
the group $G=St(1)$ of permutations which leave $1$ fixed. The matrix of
$F$ is
\begin{equation}
\left[
\begin{matrix}
t     & \beta  & \beta  & .       & .      & \beta \cr

\beta & \alpha & \delta & .       & .      & \delta \cr

\beta & \delta  & \alpha & \delta & .      & \delta \cr

\dots \cr

\beta & \delta  & .      & \delta & \alpha & \delta \cr

\beta & \delta  & .      & .      & \delta & \alpha

\end{matrix}
\right]
\end{equation}
Consider the following vectors: $v_1 = [1^s,0^{d-s}] - (s-1)\frac{\bf
j}{n} \in V^d_{k,s}$, $v_2 = [1^{s+1},0^{d-s-1}] - s\frac{\bf j}{n} \in
V^d_{k,s+1}$. The conditions $F(v_1)=0$, $F(v_2)=0$ yield

\begin{multline}
t + 2(s - 1)\beta  + (s - 1)\alpha  + (s - 1)(s -2)\delta - \cr  \phantom{.} \\
\frac {2(s - 1)(t + (s + d - 2)\beta  + (s - 1)\alpha  +(s - 1)(d - 2)\delta )}{n}   \phantom{.} + \cr \\
  \frac{(s - 1)^{2}(t + 2(d- 1)\beta  + (d - 1)\alpha  + (d - 1)(d - 2)\delta )}{n^{2}} = 0,
\end{multline}

\bigskip
\bigskip

\begin{multline}
t + 2s\beta  + s\alpha  + s(s - 1)\delta -    \phantom{.}   \frac {2s(t+ (s + d - 1)\beta  + s\alpha + s(d-2)\delta )}{n} + \cr   \phantom{.} \\
\frac {s^{2} (t + 2(d - 1)\beta  + (d - 1)\alpha  + (d - 1)(d - 2)\delta
)}{n^{2}} = 0.
\end{multline}

We also have the conditions that the sums of the diagonal and the
off-diagonal elements are equal to $0$:
\begin{equation}
t + (d-1)\alpha = 0,   (d-1)(d-2)\delta + 2(d-1)\beta=0.
\end{equation}
The determinant of the system of these $4$ equations in variables
$\alpha$, $\beta$, $\delta$ and $t$ is equal to
\begin{equation}
\frac{2}{n}(d - 1)(s - d + 1)(s - d)(d - 2 - n)
\end{equation}
and it is not equal to $0$ for the specified values of $d$, $n$ and $s$.
Hence $t=0$.

Next we prove that all off-diagonal elements of $f$ are equal to $0$.
Without a limitation of generality we can assume that $f_{12}\ne 0$. Let
$F$ be the symmetrization of $f$ by the group of permutations which map
the set $\lbrace 1,2 \rbrace$ onto itself. The matrix of $F$ is
\begin{equation}
\left[
\begin{matrix}
t     & \beta  & \beta  & .       & .      & \beta \cr

\beta & \alpha & \delta & .       & .      & \delta \cr

\beta & \delta  & \alpha & \delta & .      & \delta \cr

\dots \cr

\beta & \delta  & .      & \delta & \alpha & \delta \cr

\beta & \delta  & .      & .      & \delta & \alpha

\end{matrix}
\right].
\end{equation}
where $\alpha\ne 0$. We consider the following vectors: $v_1 =
[1^{s+1},0^{d-s-1}] - s\frac{\bf j}{n} \in V^d_{k,s}$ and $v_2 =
[0^{d-s-1},1^{s+1}] - s\frac{\bf j}{n} \in V^d_{k,s+1}$. The equations
$F(v_1) = 0$, $F(v_2) = 0$ yield

\begin{multline}
4(s-1)\beta + 2\alpha + (s-1)(s-2)\delta - 2\frac{s}{n}(2\alpha + 2(s-1 +
d-2)\beta + (s-1)(d-3)\delta)   + \cr \frac {(s - 1)^{2}(2\alpha  + 4(d -
2)\beta + (d - 2)(d - 3)\delta )}{n^2} = 0,
\end{multline}

\begin{equation}
(s + 1)s\delta  - \frac{2(s+1)s(2\beta  + (d - 3)\delta )}{n}   + \frac
{s^{2}(2\alpha  + 4(d - 2)\beta  + (d - 2)(d - 3)\delta )}{n^2} = 0.
\end{equation}
We also know that the sum of off-diagonal elements of the matrix of $F$ is
equal to $0$:
\begin{equation}
\label{off_diagonal_1_equation} 2\alpha  + 4(d - 2)\beta  + (d - 2)(d -
3)\delta = 0
\end{equation}
so, with some simplifications, the previous two equations can be rewritten
as

\begin{multline}\label{off_diagonal_2_equation}
4(s-1)\beta + 2\alpha + (s-1)(s-2)\delta   - 2\frac{s}{n}(2\alpha + 2(s +
d - 3)\beta + (s-1)(d-3)\delta) = 0,\cr \delta  - \frac{2(2\beta  + (d -
3)\delta )}{n} = 0.
\end{multline}

The systems of equations (\ref{off_diagonal_1_equation}) and
(\ref{off_diagonal_2_equation}) in variables $\alpha$, $\beta$ and
$\delta$ has determinant $8(d-2-n)(-d+s+1)$ which is not equal to $0$ for
the specified parameters $d$, $n$ and $s$ which proves that the system has
only a $0$ solution. In particular, $\alpha = 0$. We have proven that all
off-diagonal elements of the matrix of form $f$ are equal to $0$.
\end{proof}

\section{Proof of Theorem 5}

Assume that $D\subset {\mathbb R}^{d+1}$, with $0\notin \aff D$. Consider
the lattice $\Lambda = < \vert D
>_{\mathbb Z}$. Let $v\in \vert D$ be any vertex. Consider polytope $D' =
v - D$. We have $\vert D' \subset\Lambda$. One can see that there is a
positive definite quadratic form $\phi$ whose ellipsoid circumscribes
$\vert D \cap\vert D'$ and does not have any other lattice points inside
or on the boundary. Note that since $D$ satisfies the perfection property,
$\phi$ is defined up to a parameter and a scale factor.

We can now ``stretch'' the ellipsoid until it touches some other lattice
point(s) $X$. If the vertices of $D$ cannot be embedded into the set of
vertices of a centrally symmetric perfect Delaunay polytope, then the
ellipsoid will indeed meet points in lattice $\Lambda$. Proving that means
proving that the cylinder which circumscribes both polytopes $D$ and $D'$
contains other points of $\Lambda$. If that is not true, project $\Lambda$
onto $\lin D'$ along $t$. The image $\Lambda_1$ is either a lattice, in
which case the ellipsoid which circumscribes $D'$ will circumscribe a
centrally symmetric perfect Delaunay polytope in $\Lambda_1$, or the
direct sum of a lattice and a linear subspace, in which case the set of
vertices of $D'$ has dimension less than $d$-a contradiction.

The resulting ellipsoid circumscribes the polytope $D''=\conv (\vert D
\cup \vert D' \cup X)$ and has no other lattice points inside or on the
boundary. Hence $D''$ satisfies the Delaunay property. It is easy to see
that it also satisfies the perfectness property.

\section{Summary} This paper has presented two equivalent definitions of
perfect  Delaunay polytope, one related to the study of edge forms in
L-type domains, the other to a generalization of perfect forms to
inhomogeneous case. Edge forms classification is one of the cornerstone
problems in geometry of positive quadratic forms; as Voronoi noticed in
1908 \comment{\cite{bib_voronoi_rech_par}}, the study of minima of
inhomogeneous forms is strongly related to the study of perfect
homogeneous quadratic forms.

We observed that each centrally perfect Delaunay polytope can be embedded
into a centrally symmetric one. Therefore centrally symmetric perfect
Delaunay polytopes play a special role. We conjectured that the growth of
the number of perfect Delaunay polytopes in dimension $d$ up to affine
equivalence is not prohibitive; M. Dutour (2004) discovered that there is
only one perfect Delaunay polytope in dimension $6$. We have given a new
construction of infinite series of perfect Delaunay polytopes in
dimensions $d\ge 6$, based on the Gosset polytope $3_{21}$.

Further study of perfect Delaunay polytopes can have two directions: one
to estimate the growth of number of types of centrally symmetric Delaunay
supertopes, the other to find higher-dimensional analogues of the $15$,
$16$ and $22$, $23$-dimensional examples.


\begin{thebibliography}{99}
\bibitem{Baran91}  E. P. Baranovskii, Partition of Euclidean spaces into $L$%
-polytopes of certain perfect lattices. (Russian) Discrete geometry and
topology (Russian). \emph{Trudy Mat. Inst. Steklov. } \textbf{196 }
(1991), 27--46. Translated in \emph{Proc. Steklov Inst. Math. }
\textbf{196, } (1992), no. 4, 29--51.


\bibitem{CSperfect}  Conway, J. H.; Sloane, N. J. A. Low-dimensional lattices. III. Perfect forms.
 \textit{Proc. Roy. Soc. London Ser. A} \textbf{418 }(1988), no. 1854, 43--80.


\bibitem{Cox}  H.~S.~M.~Coxeter, Discrete Groups Generated by Reflections.
\emph{Annals of Math.} \textbf{35} (1934), 588-621. Repr.\ in \emph{%
Kaleidoscopes: Selected Writings of H.~S.~M.~Coxeter}, F.~A.~Sherk
\emph{et al.}, eds., Wiley, New York, 1995.

\bibitem{DG93}  M. Deza,V. P. Grishukhin, Hypermetric graphs. \emph{Quart.
J. Math. Oxford Ser. (2)} \textbf{44} (1993), no. 176, 399--433.

\bibitem{DG95}  M. Deza, V. P. Grishukhin, V. Delaunay polytopes of cut
lattices. \emph{Linear Algebra Appl.} 226/228 (1995), 667--685.

\bibitem{DG96}  M. Deza, V. P. Grishukhin, Cut lattices and equiangular
lines. Discrete metric spaces (Bielefeld, 1994). \textit{European J.
Combin.} \textbf{17} (1996), no. 2-3, 143--156.

\bibitem{DG00}  M. Deza, V. P. Grishukhin, Hypermetric two-distance spaces.
Recent advances in interdisciplinary mathematics (Portland, ME, 1997). \emph{%
\ J. Combin. Inform. System Sci.} \textbf{25} (2000), no. 1-4, 89--132.

\bibitem{DGL92}  M. Deza, V. P. Grishukhin, M. Laurent, Extreme hypermetrics
and $L$-polytopes. Sets, graphs and numbers (Budapest, 1991), 157--209,
\textit{Colloq. Math. Soc. Janos Bolyai}, \textbf{60}, North-Holland,
Amsterdam, 1992.

\bibitem{DGL95}  M. Deza, V. P. Grishukhin, M. Laurent, Hypermetrics in
geometry of numbers. \emph{\ Combinatorial optimization (New Brunswick,
NJ,
1992--1993)}, 1--109, \emph{DIMACS Ser. Discrete Math. Theoret. Comput. Sci.}%
, \textbf{20}, Amer. Math. Soc., Providence, RI, 1995.

\bibitem{DL97}  M. Deza, M. Laurent, \textit{Geometry of cuts and metrics.}
Algorithms and Combinatorics, \textbf{15.} Springer-Verlag, Berlin,
(1997).

\bibitem{DV03}  M. Dutour, F. Vallentin, \textit{Some six-dimensional rigid
forms.} Preprint \textit{math.MG/0401191} on http://ArXiv.org . To appear
in the proceedings of the 2003 Voronoi conference on analytic number
theory and spatial tessellations.


\bibitem{D04}  M. Dutour, \textit{The six-dimensional Delaunay polytopes.}
 European Journal of Combinatorics, Volume 25, Issue 4 , May 2004, Pages
535-548


\bibitem{E75}  R. Erdahl, \emph{A convex set of second-order inhomogeneous
polynomials with applications to quantum mechanical many body theory, }
Mathematical Preprint \#1975-40, Queen's University, Kingston Ontario,
(1975).

\bibitem{E92}  R. Erdahl, A cone of inhomogeneous second-order polynomials,
\textit{Discrete Comput. Geom}. \textbf{8} (1992), no. 4, 387--416.

\bibitem{E00}  R. M. Erdahl, \textit{A structure theorem for Voronoi
polytopes of lattices}, A\ talk at a sectional meeting of the AMS,
Toronto, September 22-24, (2000).

\bibitem{E01}  R. M. Erdahl, \textit{On the tame facet of the perfect domain
}$E_{6}^{\ast }$, Plenary talk at the Seventieth Birthday Celebrations for
Sergei Ryshkov, Steklov Institute, Moscow, January 24-27, (2001).

\bibitem{ER01}  R. M. Erdahl, K. Rybnikov, \emph{Supertopes,} (2002)
http://faculty.uml.edu/krybnikov/PDF/Supertopes.pdf and as math.NT/0501245
on arxiv.org.

\bibitem{KZ}  A. Korkine, G. Zolotareff, Sur les formes quadratiques,
\textit{Math. Ann}. 6 (1873), 366-389.

\bibitem{Loesch} H.-F. Loesch. Zur Reduktionstheorie von Delone-Voronoi fur
matroidische quadratische Formen. Dissertation. Univ. Bochum, 1990.

\bibitem{M03}  J. Martinet, \textit{Perfect lattices in Euclidean spaces.}
Fundamental Principles of Mathematical Sciences, \textbf{327},
Springer-Verlag, Berlin, (2003).

\bibitem{R} K. Rybnikov, \emph{REU 2001 Report: Geometry of Numbers,}  (2001)
http://www.mathlab.cornell.edu/$\tilde{~}$upsilon/REU2001.pdf

\bibitem{R99}  S. S. Ryshkov, Direct geometric description of $n$%
-dimensional Voronoi parallelohedra of the second type, (Russian) \textit{%
Uspekhi Mat. Nauk} \textbf{54} (1999), no. 1 (325), 263-264; translation
in \textit{Russian Math. Surveys} \textbf{54} (1999), no. 1, 264-265.

\bibitem{R98}  S. S. Ryshkov, On the structure of a primitive
parallelohedron and Voronoi's last problem, (Russian) \textit{Uspekhi Mat.
Nauk} \textbf{53} (1998), no. 2 (320), 161-162; translation in \textit{%
Russian Math. Surveys} \textbf{53} (1998), no. 2, 403-405.

\bibitem{V08}  G. F. Voronoi, Nouvelles applications des param\`{e}ters
continus \`{a} la th\'{e}orie des formes quadratiques, Deuxi\`{e}me
memoire, \textit{J. Reine Angew. Math.}, \textbf{134} (1908), 198-287,
\textbf{136} (1909), 67-178.

\bibitem{V52}  G. F. Voronoi, \textit{Sobranie socineniiv treh tomah},
[Collected works in three volumes], \textbf{vol. 2}, (in Russian), Kiev
(1952), Introduction and notes by B. N. Delaunay.
\end{thebibliography}
\end{document}